\documentclass{amsart}

\usepackage{amsmath}
\usepackage{amssymb}
\usepackage{amsthm}
\usepackage{amsfonts}
\usepackage{gregory_pack}
\usepackage{hyperref}
\numberwithin{equation}{section}

\long\def\symbolfootnote[#1]#2{\begingroup%
\def\thefootnote{\fnsymbol{footnote}}\footnote[#1]{#2}\endgroup}

\begin{document}
\newcommand{\cut}[1]{\langle #1 \rangle}
\newcommand{\supp}{\text{supp}\,}
\title{Exponential return times in a zero-entropy process}

\author{Paulina Grzegorek}
\address{Institute of Mathematics and Computer Science, Wroclaw University of
Technology, Wy\-brze$\dot{\mathrm z}$e Wys\-pia{\'n}\-skie\-go 27, 50-370 Wroc{\l}aw, Poland}
\email{paulina.grzegorek@pwr.wroc.pl}

\author{Michal Kupsa}
\address{Institute of Information Theory and Automation, Czech Academy of Science, Prague, Czech Republic}
\email{kupsa@utia.cas.cz}

%2010
\date{June 18, 2010}

% \thanks{The first author's research was supported by MNiSW Grant N N201 387434.

\begin{abstract}
We construct a zero-entropy weakly mixing finite-valued process with the exponential limit law for return resp. hitting times. This limit law is obtained in almost every point, taking the limit along the full sequence of cylinders around the point. Till now, the exponential limit law for return resp. hitting times, taking the limit along the full sequence of cylinders, have been obtained only in positive-entropy processes satisfying some strong mixing conditions of Rossenblatt type.
\end{abstract}
\maketitle

%------------------------------------------------------------------------
%------------------------------------------------------------------------
%------------------------------------------------------------------------
\section{Introduction}\label{sec:introduction}
In the last two decades, asymptotic laws for the return and hitting time statistics in stationary processes were intensively studied. They have been investigated mainly in the context of strong mixing properties of the process and the results are of two kinds. First under some strong mixing conditions, the limit distribution of return (resp. hitting) times to shrinking cylinders is exponential. See for instance \cite{Ab01,AV08,Pi91,GS97, GK09}. 

In these cases, the authors are taking the limit in almost every point, along the full sequence of cylinders around the point. The strong mixing conditions for processes imply positive entropy.

On the other hand, there are several classes of zero-entropy processes, which do not satisfy these strong mixing conditions and possess another limit distribution for return (resp. hitting) times.
These concern some low-complexity shifts as Sturmian shifts, linearly recurrent shifts and substitutive shifts, where the limit distributions for hitting times were proved to be piecewise linear. Chaumoitre and Kupsa \cite{CK06} proved that in the class of processes derived from rank-one systems, one can actually obtain any possible limit law for return and hitting times satisfying some very weak and natural condition described in \cite{La02} and \cite{KL05}. In particular, the exponential law can be obtained as the limit law for return and hitting times in a rank-one process. Let us recall, that rank-one processes have entropy equal to zero. 

However, most of the results from the previous paragraph are much weaker than results about the exponential limit law, since the limit laws are attained taking the limit along particular subsequences of cylinders. Although the limit law is again attained in almost every point (i.e. the limit law is still ``global''), we pick a particular sequence of cylinders around the point. This allows for coexistence of different limit laws in one process. Indeed, an example of a rank-one system, where all possible laws are realised as the limit laws for hitting times along suitably chosen subsequences of cylinders, is provided in \cite{CK06}. There are only a couple of examples of zero-entropy processes, where the limit law is attained along the full sequence of cylinders. These are the process derived from the Fibonacci shift and processes derived from two-column rank-one constructions (\cite{CK05}). In these cases the limit law for the hitting time is piecewise linear. Taking the limit along suitably chosen subsequences of cylinders also allows one to obtain a non-exponential limit distribution for positive-entropy processes, see \cite{DL}.

The question which we answer in this paper is whether there exists a zero-entropy process where the exponential law is the limit law for return and hitting times attained taking the limit along the full sequences of cylinders. Our answer is positive and we therefore demonstrate that this behavior of return and hitting times implies neither strong mixing conditions, nor positive entropy. 

Our paper is structured as follows. Section \ref{sec:preliminaries} contains our main theorem and definitions and notations needed to state it. In Section \ref{sec:definition-process}, we define the zero-entropy stationary process. Sections \ref{sec:language-analysis}, \ref{sec:close-return-hitting-times} and \ref{sec:expon-limit-distr} are devoted to the careful analysis of the process and its non-stationary measure, which is used in the construction of the process. These sections provide all technical steps needed to prove the that the process possesses the exponential limit distribution, see Corollary \ref{exp-limit-for-all-beginnings}. Theorem \ref{main-theorem} comes immediately from this corollary. The appendix contains some basic analysis of the non-stationary measure, concerning the notion of the dependency structure.

%------------------------------------------------------------------------
%------------------------------------------------------------------------
%------------------------------------------------------------------------
\section*{Acknowledgments}
We would like to thank Tomasz Downarowicz. It was his idea to construct this kind of processes to obtain exponential behavior in the class of zero-entropy processes. He encouraged us and discussions with him helped us to overcome many technicalities.

We are also grateful to El Houcein El Abdalaoui for sharing his knowledge on rank-one systems. 

The first author's research was supported by MNiSW Grant N N201 394537. The second author's research was supported by GA ASCR under the grant KJB100750901.

%------------------------------------------------------------------------
%------------------------------------------------------------------------
%------------------------------------------------------------------------
\section{Preliminaries and the main theorem}\label{sec:preliminaries}
Let $A$ be a finite set, called an alphabet. Let $A^\NN$ be the space of all sequences $x=x_0x_1\ldots$, $x_i\in A$, equipped with the $\sigma$-field $\bc$ generated by the following sets
$$[x\cut{n}]=\{y\in A^\NN:\; y_0=x_0, y_1=x_1,\ldots,  y_{n-1}=x_{n-1}\},\qquad x\in A^\NN, n\in\NN.$$
On this measurable space we consider the classical shift mapping $T:A^\NN\mapp A^\NN$
$$(Tx)_i=x_{i+1},\qquad x\in A^\NN, i\in\NN.$$
A quadruple $(A^\NN,\bc,T,\mu)$, where $\mu$ is a $T$-invariant probability measure, is called a finite-valued stationary process. This approach to define a process is characteristic for ergodic theory. The standard notion of the process, being a sequence of random variables, can be obtained if one considers the sequence of projections from $(A^\NN,\bc,\mu)$ onto $A$. 

For a measurable set $B\in\bc$ of positive measure and a point $x\in X$, we define the hitting time of $x$ in $B$ as follows,
$$\tau_B(x)=\min\{k\geq 1:\; T^kx\in B\},$$ 
where $T^k$ is the $k$-th iteration of the transformation $T$. The function $\tau_B$ can be considered as a random variable on the probability space $(X,\bc,\mu)$, or as a random variable on the conditional space $(B,\bc|B,\mu_B)$, where $\bc|B$ is the restriction of the $\sigma$-field $\bc$ to the set $B$ and $\mu_B$ is the measure defined on $\bc|B$ by the formula $\mu_B=\mu/\mu(B)$. The former variable is called the {\it hitting time} to $B$, whereas the latter one is called the {\it return time} to $B$. Poincar\'e Theorem states the return time to $B$ is almost surely finite. Since the process is assumed to be ergodic, Kac Lemma ensures the expectation of the return time to $B$ is equal to the reciprocal of the measure $\mu(B)$. 

For any $x\in A^\NN$, one can consider the sequence of suitably rescaled return or hitting times $\mu([x\cut{n}])\tau_{[x\cut{n}]}$, $n\in\NN$. The question is, whether these sequences converge in distribution. It can be rephrased in the term of weak convergence of the corresponding distribution functions. We denote the following functions from $[0,\infty]$ to $[0,1]$,
\begin{align*}
F_{x,n}(t)&=\mu\{y\in A^\NN:\;\ \mu([x\cut{n}])\tau_{[x\cut{n}]}(y)\leq t\},\\
\tilde F_{x,n}(t)&=\mu_{[x\cut{n}]}\{y\in [x\cut{n}]:\;\ \mu([x\cut{n}])\tau_{[x\cut{n}]}\leq t\}.  
\end{align*}
We say, the distribution function $F_{x,n}$ (resp. $\tilde F_{x,n}$), $n\in\NN$, weakly converges to a distribution function $F$ if $F_{x,n}(t)$ (resp. $\tilde F_{x,n}(t)$), $n\in\NN$, converges to $F(t)$ for every point $t$ of continuity of $F$. Our main result is the following:

\begin{thm}\label{main-theorem}
There exists a finite-valued weakly mixing zero-entropy process\\
$(A^\NN,\bc,T,\mu)$ such that for almost every point $x\in A^\NN$, the sequence of rescaled hitting times  $\mu([x\cut{n}])\tau_{[x\cut{n}]}$, $n\in\NN$, converges in distribution to the exponential law with parameter 1, i.e. 
$$\lim_{n\to\infty}F_{x,n}(t)=1-e^{-t}, \qquad t\geq 0.$$     
\end{thm}

An immediate consequence of the integral equation introduced in \cite{HLV05} is that the convergence in the theorem also holds for the distribution function of return times $\tilde F_{x,n}$.
\bigskip

In the rest of the section we introduce some necessary notations. For $n\in\NN$, a word (or block) of length $n$ over the alphabet $A$ is any finite sequence $u=u_0\ldots u_{n-1}$ of elements from $A$. The set of all words of length $n$ is denoted by $A^n$, the length of $u$ is denoted by $|u|$. The length of an infinite sequence $v\in A^\NN$ is formally defined to be $+\infty$.  The set of all words of all lengths is denoted by $A^*$, i.e. $A^*=\bigcup_{n\in\NN}A^n$. The concatenation of two words $u,v\in A^*$,  denoted simply by $uv$, is a word from $A^{|u|+|v|}$ such that $(uv)_i=u_i$ if $i< |u|$, and $(uv)_i=v_{i-|u|}$ if $|u|\leq i<|u|+|v|$. The concatenation of $k$ copies of a word $u$ is denoted by $u^k$.
For $u\in A^*\cup A^\NN$, $m,n\in\NN$, $m\leq n\leq |u|$, we define the word $u\cut{m,n}\in A^{n-m}$ by
$$u\cut{m,n}_i=u_{m+i},\qquad 0\leq i<m-n.$$

The language of $u$, denoted by $\lc(u)$, is a subset of $A^*$ consisting of words $u\cut{m,n}$, where $m\leq n\leq |u|$. These words are called subwords of $u$. The subword $u\cut{0,n}$ %is called the prefix of $u$ of length $n$ it 
will be denoted in the shorter way $u\cut{n}$. For a set $S\subseteq A^*\cup A^\NN$, the language of $S$ is defined as the union of the sets $\lc(u)$, $u\in S$.

For $n\in\NN$, $u\in A^n$, we denote the cylinder given by $u$ as follows,
$$[u]=\{y\in A^\NN:\; y_0=u_0, y_1=u_1,\ldots,  y_{n-1}=u_{n-1}\}.$$
This definition is an analog to the definition of the cylinder $[x\cut{n}]$ in the previous section. We will also deal with some measurable partitions of $A^\NN$. For $I\subseteq \NN$, we define the partition $\pc(I)$ as follows: Two elements $x,y\in A^\NN$ belong to the same set from $\pc(I)$ if and only if $x_i=y_i$, for every $i\in I$.
For $m,n\in\NN$, we denote 
$$\pc(m,n)=\pc([m,n)),\qquad \pc(n)=\pc([0,n)),\qquad \pc=\pc(1),$$
where $[m,n)$ and $[0,n)$ are left-closed right-open intervals of integers. For every $n\in\NN$, the partition $\pc(n)$ consists of the cylinders $[u]$, $u\in A^n$.

From now on, we reserve the symbol $A$ for three-element alphabet $\{0,1,2\}$.
%------------------------------------------------------------------------
%------------------------------------------------------------------------
%------------------------------------------------------------------------

\section{Definition of the process}\label{sec:definition-process}
In this section we will define a stationary process by averaging a non-stationary one. Both processes will have the alphabet $A=\{0,1,2\}$. 

\subsection{Non-stationary measure}\label{sec:non-stat-measure}
First, we define an auxiliary sequence $(\omega_n)^\infty_{n=1}$ of natural numbers
$$\omega_n=\max\{k\in\NN:\;2^k\mbox{ divides }n\},\qquad n=1,2,3,\ldots$$
This sequence looks as follows
$$\omega=01020103010201040102010301020105\dots$$

Proof of the following simple fact stating some properties of the sequence $\omega$ is left to the reader.
\newtheorem{fct}{Fact}
\begin{fct}\label{omega}
Let $n,m\in\NN$, $m\ge1$ and $1\le i<2^n$. Then
\begin{gather}%\begin{split}
	\omega_{n\cdot m}=	\omega_{n}+\omega_m,\label{omegasplits}\\
%\end{split}\end{equation}
%\begin{equation}
	\omega_{2^n\cdot m+i}=\omega_{i},\label{omega2}\\
	\sum^{2^n}_{j=1}\omega_j=2^n-1.\label{omegasum}
\end{gather}
The sequence $\omega$ possesses the nonoverlapping property, i.e.  $$\omega\cut{k,k+n}=\omega\cut{l,l+n}\qquad\Longrightarrow\qquad |l-k|\geq n.$$
\end{fct}
\noindent The nonoverlapping property will ensure that return times of cylinders in a process introduced below are not too small.

Let $(b_n)^\infty_{n=1}$ be a sequence of powers of $2$\symbolfootnote[2]{We will say that a number $n$ is a power of $2$ if $n=2^k$ for some natural $k$.}, increasing fast enough to satisfy the following condition
\begin{equation}
  \label{eq:b_n}
\lim_{n\to\infty}b_{n+1}/\prod^n_{i=1}b_i=\infty.  
\end{equation}
This condition is not needed until Section \ref{sec:expon-limit-distr}, where one can find the last steps of our main result's proof, see Lemma \ref{lem:b_n-increase-fast}.

We denote by $'$ the permutation on $A=\{0,1,2\}$, which permutes $1$ and $2$. This permutation is called the negation. We extend the mapping to words over $A$ and to sequences from $A^\NN$. This extension is defined letter by letter, i.e. it commutes with the concatenation and $T$. We remark, that the symbol 0 has a special role in this mapping, it is a fixpoint.

In addition, we inductively define a sequence of numbers $(a_k)^\infty_{k=1}$ and a sequence of families of words $\ac_k\subseteq A^{a_k}$, $k\ge 0$. Put $a_0=1$ and $\ac_0=\{1,2\}$. For $k\ge 0$, $a_{k+1}=2b_{k+1}a_k+2b_{k+1}-1$ and $\ac_{k+1}$ consists of all words of the form
$$u(1)0^{\omega_1}u(2)0^{\omega_2}\ldots u(b_{k+1})0^{\omega_{b_{k+1}}}(u(1))'0^{\omega_{b_{k+1}+1}}(u(2))'0^{\omega_{b_{k+1}+2}}\ldots(u(b_{k+1}))'0^{\omega_{2b_{k+1}}},$$
where $u(i)\in\ac_k$.

We define a measure $\nu$ on $\bc$ by defining its values on the set of generators of $\bc$. Let $k\in\NN$ and $u\in A^{a_k}$, then
$$\nu([u])=
\begin{cases}
1/\#\ac_k,& \mbox{if $u\in\ac_k$,}  \\  
0, & \mbox{otherwise}.
\end{cases}
$$
Since every word from $\ac_k$ has the same number of prolongations in $\ac_{k+1}$, the function $\nu$ is additive on cylinders and the definition is correct. The support of the measure (we consider the standard topology on $A^\NN$ generated by all cylinders $[u]$, $u\in A^*$) equals
$$\text{supp }\nu=\bigcap_{k\in\NN}\bigcup_{u\in\ac_k}[u].$$
The language of the support is denoted by $\lc$, i.e. $\lc=\lc(\text{supp }\nu)=\lc\left(\bigcup_{k\in\NN}\ac_k\right)$.

\begin{lem}\label{step k+l}
For $k\in\NN$, the negation is a permutation on $\ac_k$. 

For $l\geq 1$, every word $u\in\ac_{k+l}$  can be written as follows
$$u=u(1)0^{\omega_1}u(2)0^{\omega_2}\ldots u(p(k,l))0^{\omega_{p(k,l)}}, \qquad p(k,l)=\prod^l_{i=1}2b_{k+i}, \qquad u(i)\in\ac_k.$$

In particular, the support of $\nu$ is invariant under the negation and for every $x\in\text{supp }\nu$, $k\in\NN$, 
$$x=u(1)0^{\omega_1}u(2)0^{\omega_2}\ldots, \qquad u(i)\in\ac_k.$$
 \end{lem}

\begin{proof}
The negation is injective. By the inductive construction of families $\ac_k$, $k\in\NN$, they are invariant under the negation. Hence, the negation is a permutation on $\ac_k$. It implies that the second part of the lemma holds for $l=1$, $k\in\NN$.

Now, we assume the second part of the lemma holds for some $l\geq 1$ and every $k\in\NN$. Take $u\in\ac_{k+l+1}$. By the inductive assumption, 
$$u=v(1)0^{\omega_1}v(2)0^{\omega_2}\ldots v(p(l,l+1))0^{\omega_{p(l,l+1)}},\qquad v(i)\in\ac_{k+l}.$$%, for every $i\leq p(l,l+1)$.
%Using the inductive assumption again for every $v(i)\in\ac_{k+l}$, we get that for every $i\leq p(l,l+1)$, there exists sequence of words $u(j)\in\ac_k$, $(i-1)p(k,l)<j\leq ip(k,l)$, such that
%$$v(i)=u((i-1)p(k,l)+1)0^{\omega_1}u((i-1)p(k,l)+2)0^{\omega_2}\ldots u(ip(k,l))0^{\omega_{p(k,l)}}.$$
Using the inductive assumption again, there exists sequence of words $u(j)\in\ac_k$ such that for every $i\leq p(l,l+1)$
$$v(i)=u((i-1)p(k,l)+1)0^{\omega_1}u((i-1)p(k,l)+2)0^{\omega_2}\ldots u(ip(k,l))0^{\omega_{p(k,l)}}.$$
Thus, 
$$u=u(1)0^{s_1}u(2)0^{s_2}\ldots u(p(l,l+1)p(k,l))0^{\omega_{p(l,l+1)p(k,l)}},$$
where for every $1\leq i\leq p(l,l+1)$, $1\leq j\leq p(k,l)$, the following equality holds
$$s_{(i-1)p(k,l)+j}=
\begin{cases}
\omega_j,&\text{if }j<p(k,l),\\
\omega_{p(k,l)}+\omega_i,&\text{if }j=p(k,l).\\
\end{cases}
$$
Since $p(k,l)$ is a power of $2$, \eqref{omegasplits}, \eqref{omega2} and the above yields $s_{(i-1)p(k,l)+j}=\omega_{(i-1)p(k,l)+j}$. The equality $p(k,l+1)=p(k,l)p(l,l+1)$ concludes the proof.
\end{proof}

The lemma tells us at which position we can expect the occurrence of words from $\ac_k$. 
This could be reformulated in the notion of the coefficients
$$z_k(0)=0,\qquad z_k(i)=ia_k+\sum^i_{j=1}\omega_j,\qquad i,k\in\NN.$$
By the definition of $a_k$ and by equations (\ref{omegasplits})--%, (\ref{omega2}) and 
(\ref{omegasum}), we get the following arithmetic properties of the coefficients
\begin{equation}
  \label{eq:z_split}
 z_k(2^n\cdot m+i)=z_k(2^n\cdot m)+z_k(i), \qquad k,m,n,i\in\NN, i<2^n, 
\end{equation}
\begin{equation}
  \label{eq:z_k+l}
a_k=z_k(1),\quad z_{k+l}(m)=z_k(m\cdot p(k,l)), \qquad m,k,l\in\NN, k\leq l. 
\end{equation}
These properties will be often used throughout the text.% and we omit the reference. 

The following fact can be easily proved by looking at the construction of the process.
\begin{fct}\label{support-nu}
For every $k,i\in\NN$, $x\in\text{supp }\nu$, 
$$x\cut{z_k(i),z_k(i)+a_k}\in\ac_k,\qquad x\cut{z_k(i)+a_k,z_k(i+1)}=0^{\omega_{i+1}}.$$
\end{fct}

\subsection{Extension of a rank-one system}\label{sec:extension-rank-one}
To understand better the measure defined in the preceding section, we consider the continuous projection $\pi: A^\NN\to \{0,1\}^\NN$ that works letter by letter and is defined on symbols as follows: $\pi(0)=0$ and $\pi(1)=\pi(2)=1.$ By Lemma \ref{step k+l}, given $k\geq 1$, for every word $u\in\ac_k$,
$$u=u(1)0^{\omega_1}u(2)0^{\omega_2}\ldots u(p(0,k))0^{\omega_{p(0,k)}},\qquad u(i)\in\ac_0=\{1,2\}.$$
Thus, all words from $\ac_k$ have the same $\pi$-image $w(k)=10^{\omega_1}10^{\omega_2}\ldots 10^{\omega_{p(0,k)}}$. Denote $w(0)=1$. By the inductive definition of $\ac_k$, $k\in\NN$, we get 
$$
w(k+1)=w(k)0^{\omega_1}w(k)0^{\omega_2}\ldots w(k)0^{\omega_{2b_{k+1}}},\qquad  k\in\NN.
$$
Since $w(k+1)$ is a prolongation of $w(k)$, the limit $w=\lim_{k\to\infty} w(k)$ exists. The point $w$ is generic for a rank-one (non-atomic) measure on  $\{0,1\}^\NN$ (see the symbolic definition of a rank-one system in \cite{fe97}). This measure is positive on every cylinder $[u]$, $u\in\lc(w)$.
By standard Chacon arguments, see \cite{Ch69}, one can show that the rank-one system is weakly mixing and $1/2$-rigid. Since %the $\pi$-image of $\nu$ sits on $w$ (
$\nu(\pi^{-1}\{w\})=1$, we can reformulate the facts about the rank-one measure in the following way.

\begin{fct}\label{pisets}
Let $B\in\bc$ be $\pi$-measurable. Then for every $n\in\NN$,  $T^n\nu(B)\in\{0,1\}$. Moreover, Ces\`aro averages $1/n\sum^{n-1}_{i=0}T^i\nu$ weakly converge to a probability measure on the $\sigma$-field of all $\pi$-measurable sets. The limit measure is invariant and weakly mixing with respect to $T$. 

For every $u\in\lc$,  the limit measure of the set $\pi^{-1}\pi [v]$ is positive. 
\end{fct}

We will show that the Ces\`aro averages weakly converge on the whole $\sigma$-field $\bc$. The following lemma plays an important role. Its proof is introduced in the technical Section \ref{sec:language-analysis}.

%\uwaga{ZDEFINIOWAC WCZESNIEJ $\lc$!!}
\begin{lem}\label{nuconstorzero}
For $u\in\lc$, there exists $\theta(u)>0$, such that  
$$T^n\nu([u])=\theta(u)\cdot T^n\nu(\pi^{-1}\pi[u]), \qquad n\in\NN.$$
\end{lem}

\begin{lem}\label{convergenceofnu}
Measures $\frac1n\sum^{n-1}_{i=0}T^i\nu$ weakly converge.
\end{lem}

\begin{proof}
Let $u\in\lc$, then for every $n\in\NN$,
$$\frac1n\sum^{n-1}_{i=0}T^i\nu([u])=\theta(u)\frac1n\sum^{n-1}_{i=0}T^n\nu(\pi^{-1}\pi[u]),$$   
where $\theta(u)$ is the number from Lemma \ref{nuconstorzero}. Since $\pi^{-1}(\pi[u])$ is an open-closed $\pi$-measurable set, the Ces\`aro averages on the right hand side converge. Thus, the left hand side also converges. If $u\not\in\lc$, then $T^i\nu([u])=0$ for every $i\in\NN$. Thus, the Ces\`aro averages converge. 
\end{proof}

The limit measure from Lemma \ref{convergenceofnu} will be denoted by $\mu$. Fact \ref{pisets} and Lemma \ref{nuconstorzero} imply that
$$\text{ supp }\mu=\text{ supp }\nu,\qquad \lc(\text{ supp }\mu)=\lc$$
and
\begin{equation}
  \label{eq:theta}
  \theta(u)=\mu_{\pi^{-1}\pi([u])}([u]),\qquad u\in\lc.
\end{equation}

%In the rest of this section we show that the stationary process $(A^\NN,\bc,T,\mu)$ is weakly mixing and has entropy zero. %The exponential limit distribution for the return times in this process is proved in next sections.

\begin{pro}\label{weaklymixing}
The system $(A^\NN,\bc,\mu,T)$ is weakly mixing.
\end{pro}
\noindent Proof of this Proposition is presented at the end of \hyperlink{app}{Appendix}.

\begin{lem}
       The entropy of the system $(A^\NN,\bc,T,\mu)$ equals zero.
\end{lem}
\begin{proof}
It suffices to estimate from above the number of words from the language $\lc$ of length $n$, for an infinite sequence of natural numbers. Fix $k\in\NN$, $v\in\lc$, $|v|=a_k$. Since every $x\in\text{supp }\nu$ is of the form $u(1)0^{\omega_1}u(2)0^{\omega_2}\ldots$, where $u(i)\in\NN$, the word $v$ has to be a subword of $u0^{i}\tilde u$ for some $u,\tilde u\in\ac_k$ and $i\le a_k$. Hence, the number of words from $\lc$ of length $a_k$ is bounded by $2a_k^2\cdot \left(\#\ac_k\right)^2$.  The entropy of the system
       \begin{gather*}
               h(\mu)\le\lim_{k\to\infty}\frac{\log(2a_k^2\cdot \left(\#\ac_k\right)^2)}{a_k}=\lim_{k\to\infty}\frac{2\log(\#\ac_k)}{a_k}=\lim_{k\to\infty}\frac{2\prod_{i=1}^kb_i}{\prod_{i=1}^k(2b_i)}=0.
       \end{gather*}
\end{proof}

%------------------------------------------------------------------------
%------------------------------------------------------------------------
%------------------------------------------------------------------------

\section{Language analysis}\label{sec:language-analysis}

%
%We know what is $\lc(\ac_k)$, $\lc(\ac)$, $\lc_1(\ac_k)$(not block of zeros),  $\lc_{11}(\ac_k)$ (non-empty, beginning with $1$ or $2$). 
%\begin{dfn}
%Let $u\in\lc_{11}(\ac)$. The smallest $k\in\NN$ such that $u\in\lc(\ac_{k+1})$ will be called \emph{the order of $u$}. 
%\end{dfn}
For this section, let $u$ be a word from $\lc$, such that $u$ is not a block of zeros. %Since $u\in A^*$, t
There exist natural numbers $q,m_0,m_1,\ldots,m_q$ and symbols $u(i)\in\{1,2\}$, $1\leq i\leq q$, such that
$$u=0^{m_0}u(1)0^{m_1}u(2)0^{m_2}\ldots u(q)0^{m_q}.$$
The numbers and the symbols are unique and $q\geq 1$. For the rest of the section, we define $\omega_0$ to be $0$.

Denote 
\begin{align*}
\Omega(u)&=\{i\in\NN:\; \omega_{i+j}=m_j \text{ for every }1\le j<q\text{ and }\omega_i\geq m_0,\omega_{i+q}\geq m_q\},\\  
\Xi(u)&=\{z_0(i)-m_0:\; i\in\Omega(u)\}.
\end{align*}
Since $\pi(u)=0^{m_0}10^{m_1}10^{m_2}\ldots 10^{m_q}$ and $w=0^{\omega_0}10^{\omega_1}10^{\omega_2}\ldots$, we get the following lemma.
\begin{lem}\label{occ_skel}
The image $\pi(u)$ appears in the sequence $w$ only at positions from the set $\Xi(u)$.
\end{lem}

\begin{lem}\label{gaps_omega}
If $m_0=0$, then there exist $g_0,g\in\NN$, such that $g$ is a power of $2$, $g> g_0$ and 
$$\Omega(u)=\{g_0+ng,\; n\in\NN\}.$$
More precisely, $g$ is the maximum of the union $\{2^{m_j+1},\;j=0,1,\ldots,q-1\}\cup\{2^{m_q}\}$.
\end{lem}

\begin{proof}
By the definition of the sequence $\omega$, for every $j=1,2,\ldots,q-1$,
$$\Omega_j:=\{i\in\NN:\; \omega_{i+j}=m_j\}=\{(2k+1)2^{m_j}-j,\;k\in\NN\}.$$
In addition,
\begin{align*}
\Omega_q:=\{i\in\NN:\; \omega_{i+q}\geq m_q\}=\{(k+1)2^{m_q}-q,\;k\in\NN\}.
\end{align*}
Since $m_0=0$, the condition $\omega_i\geq m_0$ holds for every $i\in\NN$ and $\Omega(u)$ is equal to the intersection of the sets $\Omega_j$, $1\leq j\leq q$. Each of these sets has constant gaps between consecutive elements, which are bigger than the value of the smallest element. Their intersection is either empty (which does not happen for $u\in\lc$) or has the same property. Gaps in the intersection are also  constant and are equal to the least common multiple of the gaps in the sets themselves. Since gaps in every set $\Omega_j$, $1\leq j\leq q$, are powers of $2$, the gap $g$ in the intersection is equal to the maximum of the gaps $2^{m_q}$ and $2^{m_j+1}$, $j=1,\ldots,q-1$. 
Hence, there exists $g_0< g$, such that $\Omega(u)=\{g_0+ng, \;n\in\NN\}$.
\end{proof}
The number $g$ from the previous lemma will be called \emph{the order of $u$}.

\begin{cor} \label{gaps_xi}
Let  $m_0=0$ and $g,g_0$ be the numbers from the previous lemma. For $n\in\NN$, denote $\xi(n)=z_0(ng+g_0)$. Then $\xi(n)=z_0(ng)+z_0(g_0)$,\\
$\Xi(u)=\{\xi(n):\;n\in\NN\}$. In addition, the following inequalities hold
\begin{gather*}
z_0(ng)\leq \xi(n)<z_0(ng)+z_0(g)-\omega_g,\\
\xi(n+1)-\xi(n)\geq z_0(g),\\
\xi(n)+|u|\leq z_0(ng)+z_0(g).
\end{gather*}
\end{cor}
\begin{proof}
First part of the corollary follows from \eqref{eq:z_split} and from the fact that $g$ is a power of $2$. We prove the inequalities. Let $n\in\NN$. %, $g,g_0$ be the numbers from the previous lemma, $m_0=0$.  
Since $0\leq g_0\leq g-1$, we get
\begin{align*}
z_0(ng)\leq \xi(n)&\leq z_0(ng)+z_0(g-1)=z_0(ng)+z_0(g)-1-\omega_g\\
&<z_0(ng)+z_0(g)-\omega_g.
\end{align*}

By the definition of $z_k(i)$ and \eqref{omegasplits}
\begin{align*}
\xi(n+1)-\xi(n)&=z_0((n+1)g)-z_0(ng)=g+\sum^g_{i=1}\omega_{ng+i}=g+\omega_{(n+1)g}+\sum^{g-1}_{i=1}\omega_i\\
&=g+\omega_{n+1}+\omega_g+\sum^{g-1}_{i=1}\omega_i\geq g+\sum^g_{i=1}\omega_i=z_0(g).
\end{align*}

Finally, we assume that $\xi(n)+|u|> z_0(ng)+z_0(g)$. The right hand side of the inequality is equal to $z_0((n+1)g)-\omega_{n+1}$ (use \eqref{omegasplits}). By the definition of the sequence $w$
$$
0^{\omega_{(n+1)g}}=w\cut{z_0((n+1)g)-\omega_{(n+1)g},z_0((n+1)g)}=w\cut{z_0(ng)+z_0(g)-\omega_g,z_0((n+1)g)}.
$$
In particular, $w\cut{z_0(ng)+z_0(g)-\omega_g,z_0(ng)+z_0(g)}=0^{\omega_g}$. By the assumption and the first inequality of the Lemma, the interval $[z_0(ng)+z_0(g)-\omega_g,z_0(ng)+z_0(g)]$ is a subset of the interval $[\xi(n),\xi(n)+|u|)$. Moreover,  $w\cut{\xi(n),\xi(n)+|u|}=\pi(u)$ and $\pi(u)=10^{m_1}10^{m_2}\ldots 10^{m_q}$. It implies, that the word $0^{\omega_g}a$ is a subword of $\pi(u)$, for some $a\in\{0,1\}$. If $a=0$, then $\pi(u)$ contains $0^{\omega_g+1}$. Thus, $m_j>\omega_g$, for some $j\leq q$. By the definition of $g$, $g\geq 2^{m_j}$ and $g=2^{\omega_g}$, it is a contradiction. If $a=1$, then $0^{\omega_g}1$ occurs in $\pi(u)$ and $m_j\geq\omega_g$, for some $1\leq j<q$. For such $j$, $g>2^{m_j}$ and this is a contradiction with a fact $g=2^{\omega_g}$.  
\end{proof}
 
Now, we can prove Lemma \ref{nuconstorzero}.
\begin{proof}[Proof of Lemma \ref{nuconstorzero}]
Let $u\in\lc$. Since $[u]$ is a subset of $\pi^{-1}\pi([u])$, then by Lemma \ref{occ_skel} for every $n\not\in\Xi(u)$, $T^n\nu(\pi^{-1}\pi([u]))=T^n\nu([u])=0$. On the other hand, $T^n\nu(\pi^{-1}\pi([u]))$ equals $1$, for every $n\in\Xi(u)$. We need to prove, there is a constant $\theta(u)>0$, such that for every $n\in\Xi(u)$, $T^n\nu([u])=\theta(u)$.

There are three cases: $u$ is a block of zeros, $u$ begins with a non-zero letter and $u$ contains a non-zero letter, but not at the very beginning.

Suppose that $u$ is a block of zeros, than $[u]$ is $\pi$-measurable and the lemma holds with $\theta(u)=1$. 

Let $u$ begin with a non-zero letter. Let $g,g_0$ be the numbers from Lemma \ref{gaps_omega} and $n\in\NN$. 
By Corollary \ref{gaps_xi}, the interval $[\xi(n),\xi(n)+|u|)$ belongs to the interval $[z_0(ng),z_0(ng)+z_0(g))$. By Lemma \ref{invariance_of_nu}, $T^{z_0(ng)}\nu$ equals $\nu$ on $\pc(z_0(g))$. Thus,
% Denote by $k$, the smallest integer, such that $p(0,k)\leq g$. Since $g$ and $p(0,k)$ are powers of two, the number $g/p(0,k)$ is integer. Denote it by $l$. Since $p(0,k+1)=2b_{k+1}p(0,k)$, $l$ divides $b_{k+1}$. Suppose $l=b_{k+1}$. By the inequality \ref{sp(k,l)}, 
% $$\xi(n)=z_0(ng)+z_0(g_0)=z_0(nb_{k+1}p(k,0))+z_0(g_0)=z_k(nb_{k+1})+z_0(g_0),$$
% $$\xi(n)+|u|\leq z_0(ng)+z_0(g)=z_k(nb_{k+1})+z_k(b_{k+1}).$$
% Since $T^{z_k(nb_{k+1})}\nu$ equals $\nu$ on $\pc(0,z_k(b_{k+1}))$, we get
$$T^{\xi(n)}\nu([u])=T^{z_0(ng)+z_0(g_0)}\nu([u])=T^{z_0(g_0)}\nu([u]).$$
The last term does not depend on $n$. We put $\theta(u)=T^{z_0(g_0)}\nu([u])$.
% Let $l<b_{k+1}$. Since $g$ is a multiplier of $p(k,0)$, for every $x\in\text{supp }\nu$, there exist $v(i)\in\ac_k$, $i\leq l$, such that 

% \begin{align*}
% x\cut{z_0(ng),z_0((n+1)g)}=&x\cut{z_k(nl),z_k((n+1)l)}=v(1)0^{ng+1}v(2)0^{ng+2}\ldots v(l)0^{ng+l}\\
% =v(1)0^{ng+1}v(2)0^{ng+2}\ldots v(l)0^{ng+l}&
% \end{align*}

Let $u$ begin with $0$, but it is not a block of zeros. Then $u=0^{m_0}u(1)0^{m_1}\ldots u(q)0^{m_q}$, for some natural numbers $q,m_0,m_1,\ldots m_q$, where $q\geq 1$ and $u(i)\in\{1,2\}$, for $i\le q$. Denote $v=u(1)0^{m_1}\ldots u(q)0^{m_q}$. The set $[0^{m_0}]$ is $\pi$-measurable and has $T^n\nu$-measure $0$ or $1$ for every $n\in\NN$, thus
%$$T^n\nu([u])=\nu(T^{-n}[0^{m_0}]\cap T^{-n-m_0}[v])=\nu(T^{-n}[0^{m_0}])\nu(T^{-n-m_0}[v]).$$ 
$$T^n\nu([u])=T^n\nu([0^{m_0}]\cap T^{-m_0}[v])=T^n\nu([0^{m_0}])T^n\nu(T^{-m_0}[v]).$$ 
We have proved in the paragraph above that $T^n\nu(T^{-m_0}[v])\in\{0,\theta(v)\}$, hence $T^n\nu(u)\in\{0,\theta(v)\}$. We put $\theta(u)=\theta(v)$.

Finally, we prove that the constant $\theta(u)$ is positive. Since $u\in\lc$, there exist $n\in\NN$ and $x\in\text{supp }\nu$, such that $x\cut{n,n+|u|}=u$. It implies that $T^n\nu([u])>0$, adding the fact that $\theta(u)$ equals $T^n\nu([u])$ finishes the proof.
\end{proof}

\section{Closeness of return and hitting times}\label{sec:close-return-hitting-times}

Fix $u\in\lc$, which does not begin with zero. Let $g$ be the order of $u$. For  $t\in\NN$ we define
$$V(1,t)=\bigcup^t_{i=1}T^{-i}[u].$$
The aim of this section is to prove that the measures $\mu_{[u]}(V(1,t))$ and $\mu(V(1,t))$ are close.

The following lemma is a corollary of the fact, that for every measurable set $A\subseteq X$, the total variation distance between $\mu$ and $\mu_{X\backslash A}$ is equal to $\mu(A)$.

\begin{lem}\label{diff_measures_1}
For every $t\in\NN$,
\begin{equation*}
\left|\mu(V(1,t))-\mu_{X\backslash[0^{\omega_g}]}(V(1,t))\right|\leq\mu([0^{\omega_g}]).   
\end{equation*} 
\end{lem}

\begin{lem}\label{diff_measures_2}
For every $t\in\NN$,
\begin{equation*}
\left|\mu_{X\backslash[0^{\omega_g}]}(V(1,t))-\mu_{\pi^{-1}\pi([u])}(V(1,t))\right|\leq 2\,\mu_{\pi^{-1}\pi([u])}\!([u]).   
\end{equation*}
\end{lem}

\begin{proof}
By the definition of $w$, the block $0^{\omega_g}$ appears in $w$ exactly at positions from the following set 
$$\NN\cap \bigcup^\infty_{n=1}[z_0(ng)-\omega_{ng},z_0(ng)-\omega_g].$$
For $n\geq 1$, let $I_n$ denote the following interval of integers
$$
	I_n=(z_0(ng)-\omega_g,z_0((n+1)g)-\omega_{(n+1)g})\cap\NN=(z_0(ng)-\omega_g,z_0(ng)+z_0(g)-\omega_g)\cap\NN.
$$
Let $I=\bigcup_{n\in\NN}I_n$. We get
\begin{align*}
I&=\left\{n\in\NN:\;T^n\nu(X\backslash[0^{\omega_g}])=1\right\}\\
\NN\backslash I&=\left\{n\in\NN:\;T^n\nu(X\backslash[0^{\omega_g}])=0\right\}.
\end{align*}

Fix $n\in\NN$.  By Corollary \ref{gaps_xi}, $\xi(n)$ belongs to $I_n$. Suppose $j\in I_n$. Symmetric difference of the intervals $$
	[j+1,j+t] \;\text{ and }\;[\xi(n)+1,\xi(n)+t]
$$
consists of two intervals $J_1$ and $J_2$ (possibly empty), which lengths satisfy
$$
	|J_1|=|J_2|=\left|j-\xi(n)\right|\le|I_n|.
$$
Since $|I_n|<z_0(g)$ and the gaps in $\Xi(u)$ are at least $z_0(g)$ (Corollary \ref{gaps_xi}), we get that each set, $\Xi(u)\cap J_1$ and $\Xi(u)\cap J_2$, consists of at most one element. Hence
\begin{align*}
|T^j\nu(&V(1,t))-T^{\xi(n)}\nu(V(1,t))|=\left|\nu\Big(T^{-j}V(1,t)\Big)-\nu\Big(T^{-\xi(n)}V(1,t)\Big)\right|\\
&=\left|\nu(V(j+1,j+t))-\nu(V(\xi(n)+1,\xi(n)+t))\right|\\
&\leq \nu(V(J_1))+\nu(V(J_2))\\
&\le \mu_{\pi^{-1}\pi([u])}([u])\cdot \Big(\#(J_1\cap\Xi(u))+\#(J_2\cap\Xi(u))\Big)\\
&\leq 2\,\mu_{\pi^{-1}\pi([u])}([u]).
\end{align*}
We have  
\begin{align*}
\mu_{X\backslash[0^{\omega_g}]}(V(1,t))&=\lim_{n\to\infty}\,\frac1{\#(I\cap[0,n))}\!\sum_{i\in I\cap[0,n)}T^i\nu(V(1,t))\\
&=\lim_{n\to\infty}\,\frac1{\# (I\cap[0,z_0(ng)-\omega_g))}\sum_{i\in I\cap[0,z_0(ng)-\omega_g)}T^i\nu(V(1,t))\\
%&=\lim_{n\to\infty}\,\frac1{n(z_0(g)-1)}\!\sum^{n-1}_{m=0}\sum_{i\in I_m}T^i\nu(V(1,t))\\
&=\lim_{n\to\infty}\,\frac1{n(z_0(g)-1)}\!\sum^n_{m=1}\sum_{i\in I_m}T^i\nu(V(1,t))
\end{align*}
and
\begin{align*}
\mu_{\pi^{-1}\pi([u])}(V(1,t))&=\lim_{n\to\infty}\,\frac1{\#(\Xi(u)\cap[0,n))}\!\sum_{i\in \Xi(u)\cap[0,n)}T^i\nu(V(1,t))\\
%=&\lim_{n\to\infty}\,\frac1{\#(\Xi(u)\cap[0,z_0(ng)-\omega_g))}\sum_{i\in \Xi(u)\cap[0,z_0(ng)-\omega_g)}T^i\nu(V(1,t))\\
&=%\lim_{n\to\infty}\,\frac1n\sum^{n-1}_{m=0}T^{\xi(m)}\nu(V(1,t))=
\lim_{n\to\infty}\,\frac1n\sum^n_{m=1}T^{\xi(m)}\nu(V(1,t)).
\end{align*}
For every $m\geq 1$, $\#I_m=z_0(g)-1$. Hence,
\begin{align*}
&\hspace{-15mm}\left|\mu_{X\backslash[0^{\omega_g}]}(V(1,t))-\mu_{\pi^{-1}\pi(u)}(V(1,t))\right| \\
&\leq\limsup_{n\to\infty}\frac1{n(z_0(g)-1)}\sum^n_{m=1}\sum_{i\in I_m}\left|T^i\nu(V(1,t))-T^{\xi(m)}\nu(V(1,t))\right|\\%{z^k_{ig+g_0}+a_k-|u(1)|}\nu(V(1,t))|\\
&\leq\limsup_{n\to\infty}\frac1{n(z_0(g)-1)}\sum^n_{m=1}\sum_{i\in I_m}2\mu_{\pi^{-1}\pi([u])}([u])\leq 2\mu_{\pi^{-1}\pi([u])}([u]). 
\end{align*}
\end{proof}

\begin{lem}\label{diff_measures_3}
Let $k$ be the biggest integer, such that $p(0,k)\leq g$. Let $t\in\NN$, such that $t\leq (b_{k+2}-1)a_{k+1}$. Then
\begin{equation*}
\left|\mu_{\pi^{-1}\pi([u])}(V(1,t))-\mu_{[u]}(V(1,t))\right|\leq \mu_{\pi^{-1}\pi([u])}([u]).   
\end{equation*}
\end{lem}

\begin{proof}
Let $k$ be the integer, such that $p(0,k)\leq g<p(0,k+1)$, $t\leq (b_{k+2}-1)a_{k+1}$. The numbers  $p(0,k)$, $g$, $p(0,k+1)$ and $b_{k+1}$ are powers of $2$, $p(0,k+1)=2b_{k+1}p(0,k)$. Hence, there exist $g',g''\in\NN$, which are also powers of $2$, such that 
$$g=g'p(0,k),\quad p(0,k+1)=g''g,\qquad g'g''=2b_{k+1},\quad g'\leq b_{k+1}.$$

Fix $m\in\NN$. Since $g'$ divides $2b_{k+1}$, then $[z_k(mg'),z_k((m+1)g'))\subset[z_{k+1}(m'),z_{k+1}(m'+1))$, where $m'$ is the integer part of $mg'/2b_{k+1}$. The former interval will be denoted by $K_1'$, the latter one by $K_1$. Denote the following intervals
$$J=[\xi(m)+1,\xi(m)+t)\cap \Xi(u),$$
$$K_2'=[z_k(mg'+b_{k+1}),z_k((m+1)g'+b_{k+1})),$$
$$K_2=[z_{k+1}(m'+1),z_{k+1}(m'+b_{k+2})).$$

We will prove the following conditions:
\begin{itemize}\label{auxiliary_conditions}
\item  $T^{-\xi(m)}[u]$ is $\pc(K'_1)$-measurable with respect to $\nu$,\label{auxiliary_condition}
\item  $V(J\backslash K_1)$ is $\pc(K_2)$-measurable with respect to $\nu$,
\item  $V(J\cap K_1\backslash K_2')$ is $\pc(K_1\backslash (K'_1\cup K'_2))$-measurable with respect to $\nu$,
\item $\nu(V(J\cap K_1\cap K'_2))\leq \theta(u)$, $\nu(T^{-\xi(m)}[u]\cap V(J\cap K_1\cap K'_2))=0.$
\end{itemize}

By Corollary \ref{gaps_xi}, $\xi(m)+|u|\leq z_0((m+1)g)=z_k((m+1)g')$. Hence, the first condition holds. In addition, $z_k((m+1)g')\leq z_{k+1}(m'+1)$ and
\begin{align*}
\xi(m)+|u|+t\leq&z_{k+1}(m'+1)+(b_{k+2}-1)a_{k+1}\\
=& (m'+1)a_{k+1}+\sum^{m'+1}_{i=1}\omega_i+(b_{k+2}-1)a_{k+1}\\
\leq &(m'+b_{k+2})a_{k+1}+\sum^{m'+b_{k+2}}_{i=1}\omega_i=z_{k+1}(m'+b_{k+2}).
\end{align*}
Therefore, $J\backslash K_1\subseteq [z_{k+1}(m'+1),z_{k+1}(m'+b_{k+2})-|u|)$ and the second condition holds. 

Since $g'$ divides $2b_{k+1}$, there is $l\in\NN$, such that $$z_k((m+l)g')=z_k((m'+1)2b_{k+1})=z_{k+1}(m'+1).$$
The set $J\cap K_1\backslash K_2'$ can be written as follows,
$$J\cap K_1\backslash K_2'=\{\xi(m+i):\; 1\leq i\leq l-1, ig'\neq b_{k+1}\}.$$
For every $0\leq i\leq l-1$, $\xi(m+i)+|u|\leq z_k((m+i+1)g')$ and $T^{-\xi(m+i)}[u]$ is $\pc(z_k((m+i)g'),z_k((m+i+1)g'))$-measurable w.r.t. $\nu$. Thus, $V(J\cap K_1\backslash K_2')$ is $\pc(M)$-measurable w.r.t. $\nu$, where 
\begin{align*}
M=&\bigcup^{l-1}_{i=1}[z_k((m+i)g'),z_k((m+i+1)g'))\backslash [z_k(mg'+b_{k+1}),z_k((m+1)g'+b_{k+1}))\\
=&[z_k((m+1)g'),z_{k+1}(m'+1))\backslash K_2'.  
\end{align*}
Since $M$ is contained in $K_1\backslash (K'_1\cup K'_2)$, the third condition holds.

It remains to prove the forth condition. Let $i\in\NN$ be such that $ig'=b_{k+1}$. If $i\geq l$, then $J\cap K_1\cap K_2'$ and the forth condition holds. Assume $i\leq l-1$. Then $J\cap K_1\cap K_2'$ contains only the number $\xi(m+i)$. Hence,
$V(J\cap K_1\cap K_2')$ equals $T^{-\xi(m+i)}[u]$. Therefore, the first part of the forth condition is true. Now, assume that the second part does not hold, i.e.
$$\nu(T^{-\xi(m)}[u] \cap T^{-\xi(m+i)}[u])>0.$$
In particular, $\nu(T^{-\xi(m)}[a] \cap T^{-\xi(m+i)}[a])>0$, where $a=u_0\in\{1,2\}$. Obviously, $z_k(mg')\geq z_{k+1}(m')=z_k(m'2b_{k+1})$. In addition,
\begin{align*}
z_k(mg'+b_{k+1})=z_k((m+i)g')<z_k((m+l)g)=z_k((m'+1)2b_{k+1}).  
\end{align*}
Thus, there is $r<b_{k+1}$ such that $mg'=m'2b_{k+1}+r$. We get the following equation,
\begin{align*}
\xi(m+i)-\xi(m)=&z_0((m+i)g)-z_0(mg)=z_k((m+i)g')-z_k(mg')\\
=&z_k(m'2b_{k+1}+r+b_{k+1})-z_k(m'2b_{k+1}+r)\\
=&z_k(r+b_{k+1})-z_k(r)=b_{k+1}a_k+\sum^{b_{k+1}}_{j=r+1}\omega_j+\sum^r_{j=1}\omega_{b_{k+1}+j}\\
=&b_{k+1}a_k+\sum^{b_{k+1}}_{j=r+1}\omega_j+\sum^r_{j=1}\omega_j=b_{k+1}a_k+\sum^{b_{k+1}}_{j=1}\omega_j\\
=&b_{k+1}a_k+b_{k+1}-1=(a_{k+1}-1)/2.
\end{align*}
Let $x\in\text{supp }\nu\cap T^{-\xi(m)}[a] \cap T^{-\xi(m+i)}[a]$. By Fact \ref{support-nu}, $x\cut{z_{k+1}(m'),z_{k+1}(m')+a_{k+1}}\in\ac_{k+1}$. Thus, the first part $x\cut{z_{k+1}(m'),z_{k+1}(m')+(a_{k+1}-1)/2}$ is the negation of the second part $x\cut{z_{k+1}(m')+(a_{k+1}-1)/2,z_{k+1}(m')+(a_{k+1}-1)}$. In particular, $x_{\xi(m)}$ is the negation of $x_{\xi(m+i)}$. This contradicts the fact that $x_{\xi(m)}=x_{\xi(m+i)}=a$ (compare with Lemma \ref{dichotomy-T-nu}). Hence, the second part of the forth condition holds too.

When all the four conditions are proved, we continue with the proof of the lemma as follows. By Lemma \ref{invariance_of_nu}, the partitions $\pc(K_1)$ and $\pc(K_2)$ are $\nu$-independent. Moreover, partitions $\pc(K_1\backslash (K'_1\cup K'_2))$ and $\pc(K'_1)$ are $\nu$-independent and both are coarser than $\pc(K_1)$. It implies, that the partitions $\pc(K_1\backslash (K'_1\cup K'_2))$, $\pc(K'_1)$ and $\pc(K_2)$ are mutually $\nu$-independent. Hence, the sets $V(J\cap K_1\backslash K_2')$, $T^{-\xi(m)}([u])$ and $V(J\backslash K_1)$ are mutually $\nu$-independent. We calculate
\begin{align*}
\nu(T^{-\xi(m)}&[u]\cap V(J))\\
=&\nu(T^{-\xi(m)}[u]\cap (V(J\cap K_1\cap K'_2)\cup V(J\cap K_1\backslash K_2')\cup V(J\backslash K_1)))\\
=&\nu(T^{-\xi(m)}[u]\cap (V(J\cap K_1\backslash K_2')\cup V(J\backslash K_1)))\\
=&\nu(T^{-\xi(m)}[u])\nu(V(J\backslash (K_1\cap K_2')))
\end{align*}
and
\begin{align*}
\nu(V(J))=&\nu(V(J\cap K_1\cap K'_2)\cup V(J\backslash (K_1\cap K_2')))\\
=&\nu(V(J\backslash (K_1\cap K_2')))+\eps,
\end{align*}
where 
$$0\leq \eps=\nu\big(V(J\cap K_1\cap K'_2)\backslash V(J\backslash (K_1\cap K_2'))\big)\leq \nu(V(J\cap K_1\cap K'_2)))\leq\theta(u).$$

Since $V(\xi(m)+1,\xi(m)+t)=V(J)$ modulo $\nu$, the following inequality holds,
\begin{align*}
\left|\frac{\nu(T^{-\xi(m)}[u]\cap V(\xi(m)+1,\xi(m)+t))}{\theta(u)}-\nu(V(\xi(m)+1,\xi(m)+t)) \right|\leq\theta(u).  
\end{align*}
The inequality above holds for every $m\in\NN$. By equation (\ref{eq:theta}), 
\begin{align*}
|&\mu_{[u]}(V(1,t))-\mu_{\pi^{-1}\pi([u])}(V(1,t))|=\left|\frac{\mu_{\pi^{-1}\pi([u])}([u]\cap V(1,t))}{\mu_{\pi^{-1}\pi([u])}([u])}-\mu_{\pi^{-1}\pi([u])}(V(1,t))\right|\\
&=\lim_{m\to\infty}\left| \frac1m\sum^{m-1}_{i=0} \frac1{\mu_{\pi^{-1}\pi([u])}([u])} T^{\xi(i)}\nu([u]\cap V(1,t))-\frac1m \sum^{m-1}_{i=0} T^{\xi(i)}\nu(V(1,t))\right|\\
&\leq \limsup_{m\to\infty}\frac1m \sum^{m-1}_{i=0}\left|\frac{\nu(T^{-\xi(i)}[u]\cap V(\xi(m)+1,\xi(m)+t))} {\theta(u)}-\nu(V(\xi(m)+1,\xi(m)+t))\right|\\
&\leq \limsup_{m\to\infty}\frac1m \sum^{m-1}_{i=0}\theta(u)=\theta(u)=\mu_{\pi^{-1}\pi([u])}([u]).
\end{align*}
\end{proof}

Combining Lemmas \ref{diff_measures_1}, \ref{diff_measures_2} and \ref{diff_measures_3} gives the following result:

\begin{lem}\label{diff_measures}
Let $k$ be the biggest integer, such that $p(0,k)\leq g$. Let $t\in\NN$ be such that $t\leq (b_{k+2}-1)a_{k+1}$. Then
\begin{equation*}
\left|\mu(V(1,t))-\mu_{[u]}(V(1,t))\right|\leq 3\mu_{\pi^{-1}\pi([u])}([u])+\mu([0^{\omega_g}]).   
\end{equation*}
\end{lem}

%%%%%%%%%%%%%%%%%%%%%%%%%%%%%%%%%%%%%%%%%%%%%%%%%%%%%%%%%%%%%%%%%%%%%%%%%%%%%%%%%%%%%%%%%%%%%%%%%%%%%%%%%%%%%%%%%%%%%%
%%%%%%%%%%%%%%%%%%%%%%%%%%%%%%%%%%%%%%%%%%%%%%%%%%%%%%%%%%%%%%%%%%%%%%%%%%%%%%%%%%%%%%%%%%%%%%%%%%%%%%%%%%%%%%%%%%%%%%
\section{Exponential limit distributions for return and hitting times}
\label{sec:expon-limit-distr}

%%%%%%%%%%%%%%%%%%%%%%%%%%%%%%%%%%%%%%%%%%%%%%%%%%%%%%%%%%%%%%%%%%%%%%%%%%%%%%%%%%%%%%%%%%%%%%%%%%%%%%%%%%%%%%%%%%%%%%
%%%%%%%%%%%%%%%%%%%%%%%%%%%%%%%%%%%%%%%%%%%%%%%%%%%%%%%%%%%%%%%%%%%%%%%%%%%%%%%%%%%%%%%%%%%%%%%%%%%%%%%%%%%%%%%%%%%%%%

Denote $X'=(\text{supp }\mu)\backslash(\bigcup^\infty_{i=0}T^{-i}{0^\infty})$, $X''=X'\backslash [0]$. Let $x\in X''$. The words $x\cut{n}$, $n\geq 1$, belong to $\lc$ and begin with a non-zero letter. We denote the order of the word $x\cut{n}$ by  $g(x,n)$ (for the definition, see Lemma \ref{gaps_omega}). The maximal integer $k$ such that $p(0,k)\leq g(x,n)$ is denoted by $k(x,n)$. Since $|x\cut{n}|\leq z_0(g(x,n))$, we get that $g(x,n)$ and $k(x,n)$ tend to infinity with $n$ increasing to infinity.

\begin{lem}\label{lem:b_n-increase-fast} For every $x\in X''$
	$$
		\lim_{n\to\infty}\mu(x\cut{n})\cdot b_{k(x,n)+2}\cdot a_{k(x,n)+1}=\infty.
	$$
\end{lem}

\begin{proof}
Let $k\in\NN$. By the construction of the sequence $w$,
$$w\cut{0,a_{k+m}}=w(k)0^{\omega_1}w(k)0^{\omega_2}\ldots w(k)0^{\omega_{p(k,k+m)}},$$
for every $m\in\NN$. The number of occurrences of $w(k)$ in the word above is $p(k,k+m)$. 
Thus, by ergodicity
\begin{align*}
\mu(\pi^{-1}[w(k)])&\geq\liminf_{m\to\infty}\frac{p(k,k+m)}{a_{k+m}}=\liminf_{m\to\infty}\frac{p(k,k+m)}{p(k,k+m)a_k+\sum^{p(k,k+m)}_{i=1}\omega_i}\\
&=\liminf_{m\to\infty}\frac{p(k,k+m)}{p(k,k+m)a_k+p(k,k+m)-1}\geq \frac1{a_k+1}.
\end{align*}

Let $x\in X''$, $n\in\NN$. We denote $k=k(x,n)$, $u=x\cut{n}$ and $\xi(i)$ relates to the block $u$.
By Corollary \ref{gaps_xi} and equations \eqref{eq:z_k+l}
$$\xi(0)+|u|\leq z_0(g(x,n))\leq z_0(p(0,k+1))=a_{k+1}.$$ 
Thus, $T^{-\xi(0)}([u])$ is $\pc(a_{k+1})$-measurable. It implies that there exists $v\in\ac_{k+1}$ such that $u$ is a subword of $v$. We get,
\begin{align*}
\mu([u])\geq\mu([v])=\mu_{\pi^{-1}\pi([v])}([v])\cdot\mu(\pi^{-1}[w(k+1)])\geq \frac1{\#\ac_{k+1}\cdot(a_{k+1}+1)}.  
\end{align*}
Thus,
\begin{align*}
\mu(x\cut{n})\cdot b_{k(x,n)+2}\cdot a_{k(x,n)+1}&\geq \frac{b_{k(x,n)+2}\cdot a_{k(x,n)+1}}{\#\ac_{k(x,n)+1}\cdot(a_{k(x,n)+1}+1)}\\&=\frac{b_{k(x,n)+2}}{\prod^{k(x,n)+1}_{i=1}b_i}\cdot\frac{a_{k(x,n)+1}}{a_{k(x,n)+1}+1}.
\end{align*}
By equation (\ref{eq:b_n}), the right hand side tends to infinity, when $n$ goes to infinity. So does the left hand side.
\end{proof}

\begin{lem}\label{xcontainssomeufromack}
Let $x\in X'$. For every $k\in\NN$, there exists $u\in\ac_k$, which is a subword of $x$.
\end{lem}
\begin{proof}
Let $x\in X'$, $k\in\NN$. Since $x$ is not eventually equal to $0$, there exist a subword of $x$ of the form $u(1)0^{s_1}u(2)0^{s_2}\ldots u(p(k,l))0^{s_{p(k,l)}}$, where $u(i)\in\{1,2\}$. This word also belongs to $\lc$, thus by Lemma \ref{step k+l} %. By Fact \ref{support-nu}, 
there exists $j\in\NN$, such that $s_i=\omega_{i+j}$, for every $i\leq p(0,k)$. Surely, there is $i\leq p(0,k)$, such that $i+j$ is a multiple of $p(0,k)$. For this $i$, $s_i\geq\omega_{p(0,k)}$. Thus $0^{\omega_{p(0,k)}}$ is a subword of $x$. Since $x$ is not eventually equal to zero, there exists $m\in\NN$ such that  $x\cut{m,m+\omega_{p(0,k)}+1}$ equals $0^{\omega_{p(0,k)}}1$ or $0^{\omega_{p(0,k)}}2$. Denote $u=x\cut{m,m+\omega_{p(0,k)}+a_k}$. Then, there exist $y\in\supp \nu$ and $m'\in\NN$, such that $y\cut{m',m'+\omega_{p(0,k)}+a_k}=u$. But this could happen only for $m'\in\Xi(u)$. Since 
$$u=0^{m_0}u(1)0^{m_1}u(2)0^{m_1}\ldots u(q)0^{m_q},\quad u(i)\in\{1,2\},\quad m_0=\omega_{p(0,k)},$$
we get 
$$\Omega(u)\subseteq \{n\in\NN:\;\omega_n\geq\omega_{m_0}=\omega_{p(0,k)}\}=\{z_0(np(0,k)):\;n\in\NN\}.$$
Hence, there exists $n\in\NN$ such that $m'=z_0(np(0,k))-\omega_{p(0,k)}=z_k(n)-\omega_{p(0,k)}$ (use \eqref{eq:z_k+l}). Then
\begin{align*}
x\cut{m+\omega_{p(0,k)},m+\omega_{p(0,k)}+a_k}=&u\cut{\omega_{p(0,k)},\omega_{p(0,k)}+a_k}\\
=&y\cut{m'+\omega_{p(0,k)},m'+\omega_{p(0,k)}+a_k}\\
=&y\cut{z_k(n),z_k(n)+a_k}
\end{align*}
By Fact \ref{support-nu} the last term belongs to $\ac_k$.
\end{proof}

\begin{lem}\label{thetaismonotone}
	For all $n\in\NN$, for each $u\in\ac_n$, $\mu_{\pi^{-1}\pi([u])}([u])=\frac1{\#\ac_n}$.\\
	For all blocks $u,v\in\lc(\supp\mu)$ such that $v$ is a subword of $u$ holds 
$$\mu_{\pi^{-1}\pi([u])}([u])\le\mu_{\pi^{-1}\pi([v])}([v]).$$
\end{lem}
\begin{proof}
	Fix $n\in\NN$ and $u\in\ac_n$. By Lemma \ref{nuconstorzero} and by the fact $u\in\ac_n$, we get $\theta(u)=\nu([u])=\frac1{\#\ac_n}$. Now, let $u$ be an arbitrary block from the language $\lc(\supp\mu)$ and let $v$ be its non-empty subword. Take the cylinders defined by the words $u$ and $v$. There exists an integer $i\in\NN$ such that $[u]\subset T^{-i}[v]$. There exist $j\in\NN$ such that
 $$\mu_{\pi^{-1}\pi([u])}([u])=T^j\nu([u])\leq T^j\nu(T^{-i}[v])=T^{j+i}\nu([v])\leq\theta(v)=\mu_{\pi^{-1}\pi([v])}([v]).$$
\end{proof}

\begin{cor}\label{to0}
For $x\in X'$,
	$$
		\lim_{n\to\infty}\mu_{\pi^{-1}\pi([x\cut{n}])}([x\cut{n}])=0.%([x\cut{n}])}([x\cut{n}])=0.
	$$
\end{cor}
\begin{proof}
	Fix $\eps>0$.
	Take $k\in \NN$ so large that $\#\ac_k>\frac1{\eps}$. By Lemma \ref{xcontainssomeufromack}, there exists $N$ such that the word $x\cut{n}$ contains some $u\in\ac_k$ as a subword. By Lemma \ref{thetaismonotone}, for every $n\geq N$, the following inequalities hold 
	$$\mu_{\pi^{-1}\pi([x\cut{n}])}([x\cut{n}])\le\mu_{\pi^{-1}\pi([x\cut{N}])}([x\cut{N}])\le\mu_{\pi^{-1}\pi([u])}([u])=\tfrac1{\#\ac_k}<\eps.$$
\end{proof}
\begin{pro}
For every $x\in X''$, $t\in\NN$, 
       $$
               \lim_{n\to\infty}\sup\{|\tilde F_{[x\cut{n}]}(s)-F_{[x\cut{n}]}(s)|, 0<s\leq t\}=0.
       $$
\end{pro}
\begin{proof}
       Fix $x\in X''$, $t>0$ and $\eps>0$. Let $N\in \NN$ be such that for every $n\geq N$, the following conditions hold:
       \begin{itemize}
               \item $t\le \mu([x\cut{n}])\cdot (b_{k(x,n)+2}-1)\cdot a_{k(x,n)+1}$, (Lemma \ref{lem:b_n-increase-fast}),
               \item $\mu_{\pi^{-1}\pi ([x\cut{n}])}([x\cut{n}])<\frac\eps 4$ (Corollary \ref{to0}),
               \item $\mu([0^{g(x,n)}])<\frac\eps4$ ($\mu$ is non-atomic).
       \end{itemize}
Let $0<s\leq t$ and $n\geq N$. We simplify the notation and denote $x\cut{n}=u$, $s'=\frac s{\mu([x\cut{n}])}=\frac s{\mu([u])}$ and $V=\bigcup^{\lfloor s' \rfloor}_{i=1}T^{-i}[u]$. We have
       \begin{gather*}
               \left|\widetilde F_{[x\cut{n}]}(s)-F_{[x\cut{n}]}(s)\right|=\left|\widetilde F_{[u]}(s'\cdot\mu([u]))-F_{[u]}(s'\cdot\mu([u]))\right|\\
           =\left|\mu_{[u]}\{y\in [u]:\;  \tau_{[u]}(y)\le s'\}-\mu\{y\in X:\;\tau_{[u]}(y)\le s'\}  \right|=\left| \mu_{[u]}(V)-\mu(V)\right|.
       \end{gather*}
But $\lfloor s' \rfloor\leq (b_{k(x,n)+2}-1)\cdot a_{k(x,n)+1}$. Hence, by Lemma \ref{diff_measures},
$$
\left|\widetilde F_{[x\cut{n}]}(s)-F_{[x\cut{n}]}(s)\right|%=\left| \mu_{[u]}(V(t))-\mu(V(t))\right|
\leq
\mu([0^{\omega_{g(x,n)}}])+3\,\mu_{\pi^{-1}\pi([u])}([u])\leq\epsilon.
$$
\end{proof}

\begin{lem} Let $B\subseteq X$, $\mu(B)>0$. Then for every $t>0$,
$$\left|F_B(t)-(1-e^{-t})\right|\leq\sup\left\{|F_B(s)-\tf_B(s)|,0<s\leq t\right\}+2\cdot \mu(B).$$
\end{lem}
\begin{proof}
 For $t\geq 0$, denote $G_B(t)=\int^t_0(1-\tf_B(s))ds$ and
$$
\eta_1(t) = \tf_B(t)-F_B(t),\qquad \eta_2(t) = F_B(t)-G_B(t),\qquad\eta(t)=\eta_1(t)+\eta_2(t).
$$
We get
$$\int^t_0(1-G_B(s)-\eta(s))ds=G_B(t),\qquad t\geq 0.$$
Since $\eta_1$ and $\eta_2$ are measurable (piecewise linear) and bounded, so is $\eta$. By standard arguments from the linear ODE theory, the elementary equation above has the unique solution $G_B$ of the following form
$$G_B(t)=1-e^{-t}-\int^t_0e^{s-t}\eta(s)\,ds,\qquad t\geq 0.$$
In addition, $\eta_2(s)<\mu(B)$ for every $s\in\RR$, see \cite{HLV05}. Thus, 
\begin{align*}
|F_B(t)-(1-e^{-t})|\leq&|\eta_2(t)|+|G_B(t)-(1-e^{-t})|\leq\mu(B)+\sup_{0<s\leq t}|\eta(s)|\\
\leq&\mu(B)+\sup_{0<s\leq t}(|\eta_1(s)|+|\eta_2(s)|)\leq 2\mu(B)+\sup_{0<s\leq t}|\eta_1(s)|.  
\end{align*}
\end{proof}

\begin{cor}\label{exp-limit-for-nonzero-beginning}
For $x\in X''$, $t>0$, $\tilde F_{[x\cut{n}]}(t)$ and $F_{[x\cut{n}]}(t)$ tend to $1-e^{-t}$, whenever $n$ tends to infinity.
\end{cor}

\begin{lem}
Let $x\in X'$. There exist $N,n_0\in\NN$ such that $x_N\neq 0$ and for every $n\geq n_0$, the sets $[x\cut{n+N}]$ and $[x\cut{N,n+N}]$ are equal modulo $\mu$.
\end{lem}
\begin{proof}
Let $x\in X'$. If $x_0\neq 0$, then the lemma is trivial. Otherwise, $x=0^Nu(1)0^{m_1}u(2)0^{m_2}\ldots$, where $u(i)\in\{1,2\}$, $N,m_i\in\NN$ and $N\geq 1$. Put $g=2^N$, $n_0=g+\sum^g_{i=1}m_i$, $u=x\cut{N+n_0}$, $v=x\cut{N,N+n_0}$. We will prove that
\begin{align}\label{eq_local}
T^j\nu([u])\geq T^{j+N}\nu([v]),\qquad j\in\NN.
\end{align}
If $j+N\not\in\Xi(v)$, the inequality is trivial. Let $j+N\in\Xi(v)$. Since $w=10^{\omega_1}10^{\omega_2}\ldots$, where $\omega_i\geq N$ only when $i$ is a multiple of $g$, we get that the word
$$\pi^{-1}\pi(u)=0^N10^{m_1}10^{m_2}\ldots10^{m_g}$$
can appear in $w$ only at positions $z_0(gn)-N$, $n\in\NN$. It implies that for every $i<g$, by \eqref{omega2} holds $m_i=\omega_{gn+i}=\omega_i<N$. On the other hand $m_g\geq N$. Hence,
$$\pi^{-1}\pi(v)=10^{\omega_1}10^{\omega_2}\ldots10^{\omega_g}.$$
The word $v$ appears in $w$ exactly at the positions from the set $\Xi(v)=\{z_0(gn),n\in\NN\}$. We get that $j+N=ng$, for some $n\geq 1$. But $\omega_{ng}\geq n$ and
$$\nu(T^{-j}[0^N])=\nu(T^{-{z_0(ng)+N}}[0^N])\geq\nu(T^{-{z_0(ng)+\omega_{ng}}}[0^{\omega_{ng}}])=1.$$
Hence,
$$\nu(T^{-j}[u])=\nu(T^{-j}[0^N]\cap T^{-j-N}[v])=\nu(T^{-j-N}[v]).$$
The inequality $\ref{eq_local}$ is proved. Now, let $n\geq n_0$. Similarly, one can prove that the following equalities hold modulo $T^{j}\nu$, for $j\in\NN$
\begin{align*}
T^{-N}[x\cut{N,N+n}]&=T^{-N}[v]\cap T^{-N-n_0}[x\cut{N+n_0,N+n}]\\
&=[u]\cap T^{-N-n_0}[x\cut{N+n_0,N+n}]=[x\cut{N+n}].
\end{align*}
Finally
\begin{align*}
\mu([x\cut{N,N+n}])&=\mu(T^{-N}[x\cut{N,N+n}])\\
&=\lim_{k\to\infty}\frac1k \sum^{k-1}_{j=0}T^j\nu(T^{-N}[x\cut{N,N+n}])\\
%&=\lim_{k\to\infty}\frac1k \sum^{k-1}_{j=0}T^{j+N}\nu([x\cut{N,N+n}])\\
&= \lim_{k\to\infty}\frac1k \sum^{k-1}_{j=0}T^j\nu([x\cut{N+n}])=\mu([x\cut{N+n}]).
\end{align*}
%Since $T^{-N}[v]$ is a subset of $[u]$, the sets are equal modulo $\mu$.

\end{proof}

\begin{cor}\label{exp-limit-for-all-beginnings}
For $x\in X'$, $t>0$, $\tilde F_{[x\cut{n}]}(t)$ and $F_{[x\cut{n}]}(t)$ tend to $1-e^{-t}$, whenever $n$ tends to infinity.
\end{cor}
\begin{proof}
Let $x\in X'$ and $N,n_0\in\NN$ be numbers, for which the previous lemma is satisfied. Denote $y=T^Nx$. Then $y\in X''$ and for every $n\geq n_0$, $[x\cut{N+n}]$ equals $[y\cut{n}]$ modulo $\mu$. Thus, for every $n\geq n_0$, the functions $F_{x,N+n}$ and $F_{y,n}$ are the same. By Corollary \ref{exp-limit-for-nonzero-beginning}, for every $t>0$, $F_{y,n}(t)$ tends to $1-e^{-t}$, and so does the sequence $F_{x,n}(t)$.% tends to  $1-e^{-t}$ too.
\end{proof}

%\appendix

\section*{Appendix. Dependency structure on coordinates}\label{sec:depend-struct-coord}\hypertarget{app}{}
%, explicit form

The mutual dependences among partitions $T^{-j}\pc$, $j\in\NN$,  with respect to the measure $\nu$,  can be described through a non-reflexive, symmetric and transitive relation on the set of coordinates $\NN$, which is an equivalence relation on $Q=\{j\in\NN:\;w_j=1\}=\{z_0(i):\;i\in\NN\}$. 

% For $k,i,j\in\NN$, let us denote the left-closed right-open intervals of integers 
% $$
% 	I_k(i)=[z_k(i),z_k(i+1))
% $$
% \uwaga{but more usefull will be:}
% \begin{equation}\label{dfnI}
% I_k(i,j)=[z_k(i),z_k(j))
% \end{equation}

First, we define the \emph{direct dependency} symmetric relation $D$. We say, that coordinates $i,j\in\NN$ are directly dependent ($(i,j)\in D$), if $i,j\in Q$ and there exist $k,m\in\NN$ such that 
$$
	i,j\in [z_k(m),z_k(m)+a_k)\quad\mbox{ and }\quad|i-j|=\tfrac12(a_{k}-1).
$$
The dependency relation $\ovd$ on $\NN$ is the transitive envelope of $D$. We remark that this relation is a subset of $Q^2$ and it is an equivalence on $Q$. For $I\subseteq \NN$,  we denote
$$
	D(I)=\bigcup_{i\in I}\{j\in\NN:\; (i,j)\in D\}, \qquad \ovd(I)=\bigcup_{i\in I}\{j\in\NN:\; (i,j)\in\ovd\}.
$$
Since $\ovd$ is an equivalence on $Q$, the family $\dc=\{\ovd(i):\;i\in\NN\}$ is a partition of $Q$. 
The number of sets from $\dc$ crossing $I$ is denoted by $d(I)$. Again, for an interval of integers $I=[m,n)$ we omit brackets and write $\ovd(m,n)=\ovd(I)$ and $d(m,n)=d(I)$. 

\begin{lem}\label{dichotomy-T-nu}
 For $j\in\NN\backslash Q$, $T^j\nu([0])=1$. In particular, the partition $T^{-j}\pc$ is trivial with respect to $\nu$. For $j\in Q$
 $$T^j\nu([1])=T^j\nu([2])=1/2.$$

If $i,j\in\NN$, $(i,j)\in D$, $a\in\{1,2\}$, then $\nu(T^{-i}[a]\cap T^{-j}[a'])=\nu(T^{-i}[a])=1/2$. In particular, 
$$T^{-i}\pc=T^{-j}\pc=\pc(\ovd(i))\ \ \ \mod \nu.$$
 \end{lem}
 \begin{proof}
The first part of the lemma follows from definitions. Since the measure $\nu$ is invariant under the negation, $\nu(T^{-i}[1])= \nu(T^{-i}[2])=1/2$ for every $i\in Q$.

Let $(i,j)\in D$. In order to prove, that $\nu(T^{-i}[a]\cap T^{-j}[a'])$ equals $\nu(T^{-i}[a])$, $a=1,2$, we prove that for every $x\in\text{supp } \nu$, $x'_i=x_j$. Let $k,m\in\NN$ be such that $i,j\in [z_k(m),z_k(m)+a_k)$ and $|i-j|=\tfrac12(a_{k}-1)$. By Fact \ref{support-nu}, $x\cut{z_k(m),z_k(m)+a_k}=u\in\ac_k$. Denote $i'=i-z_k(m)$ and $j'=j-z_k(m)$. Since $u\cut{\tfrac12(a_{k}-1),a_{k}-1}$ is the negation of $u\cut{\tfrac12(a_{k}-1)}$, we get
$x_i'=u_{i'}'=u_{j'}=x_j$.

It also implies that  $T^{-i}\pc$ equals $T^{-j}\pc$ modulo $\nu$. We consider the relation on $\NN$ defined as follows: $(i,j)\in \NN^2$ is in the relation  if $T^{-i}\pc$ equals $T^{-j}\pc$ modulo $\nu$. Since this relation is transitive and contains $D$, it contains $\ovd$ too. Thus, for every $(i,j)\in\ovd$, $T^{-i}\pc$ equals $T^{-j}\pc$ modulo $\nu$. It implies, that  $\pc(\ovd(i))$ equals $T^{-i}\pc$ (mod $\nu$).
\end{proof}

We denote the following sets, $k\in\NN$,
\begin{gather*}
\Gamma(k)=\{\mathbb{\gamma}\in\NN^{\NN}:\; (\forall i<k)\; \gamma_i<2b_{i+1}, (\forall i\geq k)\; \gamma_i=0\},\\
\Gamma'(k)=\{\mathbb{\gamma}\in\Gamma(k):\; (\forall i\in\NN)\; \gamma_i<b_{i+1}\},\\
\Gamma=\bigcup_{k\in\NN}\Gamma(k), \qquad\Gamma'=\bigcup_{k\in\NN}\Gamma'(k).
\end{gather*}
%Denote $\Gamma=\bigcup_{k\in\NN}\Gamma(k)$, $\Gamma'=\bigcup_{k\in\NN}\Gamma'(k)$. 
We define the mapping $\phi:\Gamma\mapp Q$ by 
$$\phi(\gamma)=\sum^{\infty}_{i=0}z_i(\gamma_i).$$ 
The mapping is well defined what follows by two facts. The first is $z_i(0)=0$, for every $i\in\NN$ and the second is the equality which can be deduced from equations (\ref{eq:z_split}) and (\ref{eq:z_k+l})
$$
	\sum_{i=0}^\infty z_i(\gamma_i)=z_0\left(\sum_{i=0}^\infty \gamma_i\cdot p(0,i)\right),\qquad\gamma\in\Gamma.
$$
By the same arguments, the mapping $\phi$ is a bijection. Moreover, for every $k\in\NN$, $\phi$ maps bijectively $\Gamma(k)$ onto $Q\cap [0,a_k)$. The following lemma is an easy observation.

\begin{lem}\label{same-z_k}
Let $k\in\NN$, $\gamma,\gamma'\in\Gamma$. The following conditions are equivalent
\begin{itemize}
\item $\gamma_i=\gamma_i'$ for every $i\geq k$,
\item there exists $m\in\NN$ such that $\phi(\gamma)$ and $\phi(\gamma')$ belong to $[z_k(m),z_k(m+1))$. 
\end{itemize}
\end{lem}

\begin{lem}
Let $k\in\NN$, $\gamma,\gamma'\in\Gamma$. Assume that $\gamma_i=\gamma'_i$, for every $i\neq k$, $\gamma'_k-\gamma_k=b_{k+1}$. Then $(\phi(\gamma),\phi(\gamma'))\in D$.
\end{lem}

\begin{proof}
Let the number $k\in\NN$ and the sequences $\gamma,\gamma'\in\Gamma(k)$ satisfy the assumptions of the lemma. Denote by $m$ the integer which satisfies the condition $z_{k+1}(m)=\sum^{\infty}_{i=k+1}z_i(\gamma_i)$.  Then the differences $\phi(\gamma)-z_{k+1}(m)$ and $\phi(\gamma')-z_{k+1}(m)$ are bounded from above by the term
$$\sum^{k}_{i=0}z_i(2b_{i+1}-1)<a_{k+1}.$$
Thus, $\phi(\gamma),\phi(\gamma')\in[z_{k+1}(m),z_{k+1}(m+1))$. In addition, by \eqref{eq:z_split} (use the fact that $\gamma_k<b_{k+1}$) and by \eqref{omegasum}
\begin{align*}
\phi(\gamma')-\phi(\gamma)&=z_k(\gamma_k+b_{k+1})-z_k(\gamma_k)=z_k(\gamma_k)+z_k(b_{k+1})-z_k(\gamma_k)\\
&=z_k(b_{k+1})=b_{k+1}a_k+b_{k+1}-1=(a_{k+1}-1)/2.
\end{align*}
Thus, $(\phi(\gamma),\phi(\gamma'))\in D$.
\end{proof}
Using the fact, that $\ovd$ is a transitive envelope of $D$, one gets the following lemma.

\begin{lem}
Let $\gamma,\gamma'\in\Gamma$ and $\gamma_i'-\gamma_i\in\{0,b_{i+1}\}$, for every $i\in\NN$. Then $(\phi(\gamma),\phi(\gamma'))\in \ovd$.
\end{lem}

This lemma has the following corollary.
\begin{cor}
  For every $\gamma\in\Gamma'$
$$D(\phi(\gamma))=\left\{\sum^{\infty}_{i=0}z_i(\gamma_i+\alpha_ib_{i+1}):\; \alpha_i\in\{0,1\}, \sum_{i\in\NN}\alpha_i<\infty\right\}.$$
The sets $D(\phi(\gamma))$, $\gamma\in\Gamma'$, are pairwise disjoint and cover all $Q$.
\end{cor}
Combining this corollary with Lemma \ref{same-z_k} gives the following technical lemmas.
\begin{lem}\label{gap_for_dep}
Let $k,m,m',j,j'\in\NN$ be such that $(j,j')\in\ovd$, $j\in [z_k(m),z_k(m+1))$ and $j'\in [z_k(m'),z_k(m'+1))$. Then $b_{k+1}$ divides $|m-m'|$. 
Moreover
\begin{align*}
d(z_k&(m),z_k(m+b_{k+1}))=\sum^{b_{k+1}-1}_{i=0}d(z_k(m+i),z_k(m+i+1)),\\
d(z_k(m+&1),z_k(m+2b_{k+1}))\\
=&d\big([z_k(m+1),z_k(m+b_{k+1}))\cup[z_k(m+b_{k+1}+1),z_k(m+2b_{k+1}))\big)\\
&+d([z_k(m+b_{k+1}),z_k(m+b_{k+1}+1))).
\end{align*}
\end{lem}
\begin{lem}\label{dependency_in_interval}
For $k,m\in\NN$, $\gamma\in\Gamma'(k)$,  
$$D(z_k(m)+\phi(\gamma))\cap[z_k(m),z_k(m+1))=\left\{z_k(m)+\sum^{k-1}_{i=0}z_i(\gamma_i+\alpha_ib_{i+1}),\alpha_i\in\{0,1\}\right\}.$$
The sets $D(z_k(m)+\phi(\gamma))$, $\gamma\in\Gamma'(k)$, are pairwise disjoint and cover all $Q\cap [z_k(m),z_k(m+1))$. 

In particular, $d(z_k(m),z_k(m+1))=p(0,k)/2^k.$
\end{lem}

\begin{lem}
Partitions $\pc(E)$, where $E$ runs over all equivalence classes of $\ovd$, are mutually $\nu$-independent.
\end{lem} 
\begin{proof}
Let $k\in\NN$. Elements of the partition $\pc(a_k)$ consists of the cylinders given by the words from $\ac_k$. By the definition of $\nu$ all elements of partition $\pc(a_k)$ have the same measure. Thus, we can easily calculate the entropy of the partition
$$
	H(\pc(a_k),\nu)=\log_2(\#\ac_k)=\log_2\left(2^{\prod^k_{i=1}b_i}\right)=\prod^k_{i=1}b_k.
$$
On the other hand, the sets $D(j)$, $j\in\phi(\Gamma'(k))$, cover $Q\cap [0,a_k)$. Hence, the joining of the partitions $\pc(D(j))$,  $j\in\phi(\Gamma'(k))$, is finer than $\pc(a_k)$. In addition, by Lemma \ref{dichotomy-T-nu}, $H(\pc(D(j)),\nu)=1$, for every $j\in\phi(\Gamma'(k))$. Hence,
\begin{align*}
H(\pc(a_k),\nu)\leq& H\left(\bigvee_{j\in\phi(\Gamma'(k))}\pc(D(j)),\nu\right)\leq\sum_{j\in\phi(\Gamma'(k))}H(\pc(D(j)),\nu)\\
=&\ \#\phi(\Gamma'(k))=\prod^k_{i=1}b_k.
\end{align*}
Since the first term of the inequality above equals the last one, all the terms above are equal. In particular, the partitions $\pc(D(j))$,  $j\in\phi(\Gamma'(k))$, are mutually $\nu$-independent and this is true for every $k\in\NN$. Since $D(j)$, $j\in\bigcup_{k\in\NN}\Gamma'(k)$, are all classes of the equivalence $\ovd$, the lemma holds.
\end{proof}
\begin{cor}\label{entropy=d}
For every $I\subseteq\NN$, $H(\pc(I),\nu)=d(I)$.
\end{cor}

\begin{lem}\label{invariance_of_nu}
For every $n,m\in\NN$, the measure $T^{z_0(m2^n)}\nu$ is equal to $\nu$ on $\pc(z_0(2^n))$. 
In particular, for every $k,m\in\NN$,  $T^{z_k(m)}\nu$ is equal to $\nu$ on $\pc(a_k)$. 

Moreover, the partitions
$\pc(z_k(m),z_k(m+1))$ and
$$\pc\big([z_k(m+1),z_k(m+b_{k+1}))\cup[z_k(m+b_{k+1}+1),z_k(m+2b_{k+1}))\big)$$
are $\nu$-independent.
\end{lem}
\begin{proof}
Let $m,n\in\NN$. Let $k$ denote the biggest integer such that $p(0,k)\leq 2^n$. Then $2^n=lp(0,k)$ for some $l\in\NN$, such that $l$ divides $b_{k+1}$. Of course, $l$ is a power of $2$. Let $x\in\text{supp }\nu$. Using \eqref{omegasplits} we get
\begin{align*}
x\cut{z_0(m2^n),z_0((m+1)2^n)}=&x\cut{z_k(ml),z_k((m+1)l)}\\
=&u(1)0^{\omega_{ml+1}}u(2)0^{\omega_{ml+2}}\ldots u(l)0^{\omega_{(m+1)l}}\\
=&u(1)0^{\omega_1}u(2)0^{\omega_2}\ldots u(l)0^{\omega_l}0^{\omega_{(m+1)}}
\end{align*}
for some $u(1),u(2),\ldots,u(l)\in\ac_k$. 

Denote
$$U=\{u(1)0^{\omega_1}u(2)0^{\omega_2}\ldots u(l)0^{\omega_l}:\; u(i)\in\ac_k\}.$$ For every $v\in U$, $|v|=z_k(l)=z_0(2^n)$. Thus,  
$$\pc(z_0(m2^n),z_0(m2^n)+z_0(2^n))=\{T^{-m2^n}[v], v\in U\} \text{ modulo } \nu .$$
This partition is equal to $\pc(z_k(ml),z_k((m+1)l))$ modulo $\nu$, because the partition $\pc(z_k(ml)+z_k(l),z_k((m+1)l))$ is trivial.
By Lemma \ref{entropy=d}, Lemma \ref{gap_for_dep} and Lemma \ref{dependency_in_interval},
\begin{align*}
\log_2(\# U)&=\log_2((\#\ac_k)^l)=\log_2((2^{\prod^k_{i=1}b_k})^l)=lp(0,k)/2^k\\
&=\sum^{l-1}_{i=0}H(\pc(z_k(ml+i),z_k((ml+i+1)l)),\nu)\\
&=H(\pc(z_k(ml),z_k((m+1)l)),\nu).
\end{align*}
This is possible only when $\nu$ is constant on the partition $\{T^{-m2^n}[v], v\in U\}$. Hence, for $u\in U$, $\nu(T^{-m2^n}[v])=1/\#U$. The value $1/\#U$ does not depend on $m$, therefore, the first part of lemma holds.

Now, let $k,m\in\NN$. By Corollary \ref{entropy=d} and Lemma \ref{gap_for_dep}, 
the partitions
$\pc(z_k(m),z_k(m+1))$ and
$$\pc([z_k(m+1),z_k(m+b_{k+1}))\cup[z_k(m+b_{k+1}+1),z_k(m+2b_{k+1})))$$
are $\nu$-independent.
\end{proof}

\begin{proof}[Proof of Proposition \ref{weaklymixing}]
To avoid lengthy expressions, we denote 
$$\alpha(n,B,C)=\frac1n\sum^n_{j=1} \left| \mu(B\cap T^{-j}C)-\mu(B)\mu(C) \right|,\qquad n\in\NN,\;\, B,C\in\bc.$$
Let $u,v\in\lc(\supp\mu)$, $\eps>0$. We will prove that there exists $n_1\in\NN$ such that for every $n\geq n_1$, $\alpha(n,[u],[v])<5\eps.$

Let $B=\pi^{-1}\pi([u])$ and $C=\pi^{-1}\pi([v])$. By a standard Chacon argument one can show that the rank-one factor is weakly mixing. Thus, there exists $n_0$ such that $\alpha(n,B,C)<\eps$ for every $n\geq n_0$.

Let
$$
E=\{(i,j)\in\NN^2:\; \ovd([i,i+|u|))\cap[i+j,i+j+|v|)\neq\emptyset\}.
$$
A pair $(i,j)$ belongs to $E$ if and only if $\pc[i,{i+|u|})$ and $\pc[i_j,{i+j+|v|})$ are $\nu$-dependent. In particular, for $(i,j)\not\in E$ holds
$$\nu(T^{-i}[u]\cap T^{-i-j}[v])=\nu(T^{-i}[u])\cdot\nu(T^{-i-j}[v]).$$
In addition, we denote
%\begin{center}
%$\begin{array}{r@{}c@{}l@{}r@{}cl}
%E_i&=&\{j\in\NN:\; (i,j)\in E\},\qquad E_{i,n}&=&E_i\cap [1,n],\\
%E_j'&=&\{i\in\NN:\; (i,j)\in E\},\qquad E_{j,m}'&=&E_j\cap [1,m].
%\end{array}$
%\end{center}\vspace{-7mm}
$$
	E_{i,n}=\left\{j\le n:\; (i,j)\in E\},\qquad E_{j,m}'=\{i\le m:\; (i,j)\in E\right\}.
$$

By Lemma \ref{dependency_in_interval}, the upper Banach density of $\ovd(i)$ is bounded by $2^k/a_k$, for every $k\in\NN$. Hence $\ovd(i)$ has Banach density equal to zero. Hence, the sequence $\#E_{i,n}/n$ converges to $0$ uniformly in $i$, when $n$ goes to infinity. Therefore, there is $n_1\geq n_0$ such that for every $n\geq n_1$ and every $i\in\NN$, densities $\#E_{i,n}/n$ are less than $\eps^2$. 

Let $n\geq n_1$. For every $m\in\NN$ 
$$
	\bigcup_{i=1}^m\left(\{i\}\times E_{i,n}\right)=E\cap \left([1,m]\times [1,n]\right).
$$
%the columns $\{i\}\times E_{i,n}$, $1\leq i\leq m$, sup up to the set $E\cap [1,m]\times [1,n]$.
Hence, the density of $E$ in the rectangle $[1,m]\times [1,n]$ is less than $\eps^2$. Thus, there are less than $\eps\cdot n$ rows among $[1,m]\times \{j\}$, $1\leq j\leq n$, such that the density of $E\cap([1,m]\times \{j\})$ is bigger or equal to $\epsilon$. In other words, the set
$$M_m=\{j\leq n:\; \#E_{j,m}'/m\geq\eps\}$$ 
satisfies $\#M_m/n<\eps$.

From now on, we will use the notation $a=b+o(e)$ instead of $|a-b|\leq\eps$. 
Take $m\in\NN$ such that for every $j\leq n$ holds
\begin{align*}
\frac1m\sum^m_{i=1}T^i\nu(B\cap T^{-j}C)&=\mu(B\cap T^{-j}C)+o(\eps),\\
\frac1m\sum^m_{i=1}T^i\nu([u]\cap T^{-j}[v])&=\mu([u]\cap T^{-j}[v])+o(\eps).
\end{align*}

If $j\in [1,n]\backslash M_m$, $i\in [1,m]\backslash E'_{j,m}$, then
\begin{align*}
T^i\nu([u]\cap T^{-j}[v])&= \nu(T^{-i}[u]\cap T^{-i-j}[v])=\nu(T^{-i}[u])\nu(T^{-i-j}[v])\\
&=\mu_B([u])\nu(T^{-i}B)\mu_C([v])\nu(T^{-i-j}C)\\
&=\mu_B([u])\mu_C([v])\nu(T^{-i}B\cap T^{-i-j}C)\\
&=\mu_B([u])\mu_C([v])T^i\nu(B\cap T^{-j}C).
\end{align*}
Hence, for $j\in [1,n]\backslash M_m$,
\begin{align*}
\mu([u]\cap T^{-j}[v])&=o(\eps)+\frac1m\sum^m_{i=1}T^i\nu([u]\cap T^{-j}([v]))\\
&=o(\eps)+\mu_B([u])\mu_C([v]) \frac1m\sum^m_{i=1}T^i\nu(B\cap T^{-j}C)\\
&\quad +\frac1m\sum_{i\in E'_{j,m}}\left(T^i\nu([u]\cap T^{-j}[v])-\mu_B([u])\mu_C([v])T^i\nu(B\cap T^{-j}C)\right)\\
&=o(2\eps)+\mu_B([u])\mu_C([v]) \frac1m\sum^m_{i=1}T^i\nu(B\cap T^{-j}C)\\
&=o(3\eps)+\mu_B([u])\mu_C([v])\mu(B\cap T^{-j}C).
\end{align*} 
Finally, for every $n\geq n_1$,
\begin{align*}
\alpha(n,[u],[v])\leq& 
\frac1n\sum^n_{j=1}\left| \mu([u]\cap T^{-j}[v])-\mu_B([u])\mu_C([v])\mu(B\cap T^{-j}C)\right|\\
& +\mu_B([u])\mu_C([v])\frac1n\sum^n_{j=1}\left|\mu(B\cap T^{-j}C)-\mu(B)\mu(C)\right|\\
\leq& \mu_B([u])\mu_C([v])\alpha(n,B,C)+3\eps\\
&+ \frac1n\sum_{j\in M_m}\left| \mu([u]\cap T^{-j}[v])-\mu_B([u])\mu_C([v])\mu(B\cap T^{-j}C)\right|\\
\leq& 4\eps+\#M_m/n< 5\eps.
\end{align*}
\end{proof}

\edb